\documentclass[11pt]{amsart}   %%%%%%%Please do not change
\newtheorem{theorem}{Theorem}[section]
\newtheorem{lemma}[theorem]{Lemma}

\newtheorem{corollary}[theorem]{Corollary}
\newtheorem{example}[theorem]{Example}
\newtheorem{question}[theorem]{Question}
\theoremstyle{definition}

\theoremstyle{remark}

\numberwithin{equation}{section}

\begin{document}

\title[Topology Proceedings  Example Article]%
{Distributive lattice orderings \\ and Priestley Duality}

%    Information for first author:
\author{Michel Krebs}
\address{Mathematical Institute, University of Bern, Sidlerstrasse 5,
3012 Bern, Switzerland}
%    Current address (if needed): 
%\curraddr{}
\email{michel.krebs@math.unibe.ch}

%    Information for second author (if needed): 
\author{Dominic van der Zypen}
\address{Allianz Suisse Insurance Company, 3001 Bern, Switzerland}
\email{dominic.zypen@gmail.com}

%    General info


\begin{abstract} 
The ordering relation of a bounded distributive lattice $L$ is a
(distributive) $(0,1)$-sublattice of $L\times L$. 
This construction gives rise to
a functor $\Phi$ from the category of bounded distributive lattices to
itself. We examine the interaction of $\Phi$ with Priestley duality and
characterise those bounded distributive lattices $L$ such that there is
$K$ with $\Phi(K) \cong L$.
\end{abstract}

\maketitle

%------------------------ now paper starts ---------------------

%%%%%%%%%%%%%%%%%%%%% new commands
\def\Ip{{\mathcal I}_p}
\def\D{{\mathcal D}_{01}}
\def\X{{\mathcal X}}
\def\E{{\mathcal E}}
\def\P{{\mathcal P}}
%%%%%%%%%%%%%%%%%%%%% end of commands

\section{Some conventions and definitions}
For any poset $P$ we say that $A\subseteq P$ is a {\em lower set} or 
{\em down-set}
if $a\in A$, $x\in P, x\leq a$ imply $x\in A$. The dual notion is that
of an {\em up-set}. 

We assume all lattices to be distributive and bounded by 0,1 such that $0\neq 1$.
A nonempty down-set $I$ of a bounded 
distributive lattice $L$ is said
to be an {\em ideal} if $a,b \in I$ implies $a\vee b\in I$. 
An up-set with the dual property is called a {\em filter}.
Moreover, $I$
is a {\em prime ideal} if $I\neq L$ and if $a,b \in L\setminus I$ 
implies $a\wedge b\in L\setminus I$. Note that a down-set of $L$ is a prime ideal if
and only if its complement is a filter.

Let $L$ be a bounded distributive lattice. Then by $\Ip(L)$ we denote the set of 
all prime ideals of $L$. Suppose that $a\in L$, then we define $$ X^L_a =
\{I\in \Ip(L) : a\notin I\}.$$ When no confusion arises, we omit the 
superscript and write $X_a$.
%------------------------------------
\section{Priestley duality}
In \cite{Pr70}, Priestley proved that the category $\D$
of bounded distributive lattices with $(0,1)$-preserving lattice
homomorphisms and the category $\P$ of compact totally
order-disconnected spaces (henceforth referred to as {\it
Priestley spaces}) with order-preserving continuous maps are
dually equivalent. A compact {\it totally order-discon\-nec\-ted
space} $(X;\tau ,\leq )$ is a poset $(X;\leq )$ endowed with a
compact topology $\tau$ such that, for $x$, $y\in X$, whenever
$x\not\geq y$, then there exists a clopen decreasing set $U$ such
that $x\in U$ and $y\not\in U$. We usually refer to a Priestley space
by its ground set only when there is no ambiguity about the topology
and the ordering relation being used. In the following we briefly describe the
pair of contravariant functors connecting $\D$ and $\P$.

The functor $\X:{\D}\to
{\P}$ assigns to each object $L$ of ${\D}$ a
Priestley space $(\Ip(L);\tau (L),\subseteq)$, where $\Ip(L)$ is the
set of all prime ideals of $L$ and the topology $\tau (L)$ is
given by the following subbasis of $\Ip(L)$:
$$\{X_a: a\in L\} \cup \{\Ip(L)\setminus X_a: a\in L\}.$$

As lined out in \cite{DaPr02}, this topology is compact and totally
order-disconnected; moreover it turns out that the collection of clopen
down-sets consists exactly of the $X_a$ (and the collection of clopen
up-sets are their complements).

For a $(0,1)$-preserving lattice homomorphism $f:L\to K$ we define
$\X(f):\X(K) \to \X(L)$ by $\X(f)(I) = f^{-1}(I)$ for all $I\in \X(K)=
\Ip(K)$, using the fact that preimages of prime ideals are prime ideals.

The functor $\E:{\mathcal{P}}\to {\mathcal{D}}$ assigns to each
Priestley space set $\E(X)$
of all clopen down-sets of $X$ ordered by set inclusion, which gives
rise to a bounded distributive lattice.

On the level of morphisms, i.e. order-preserving continuous maps,
$\E$ again works with preimages.

An excellent introduction to Priestley duality can be found in \cite{DaPr02}.
%------------------------------------
\section{The functor $\Phi$}
Let $(P,\leq_P)$ be a nonempty poset. The cartesian product $P\times P$ of the 
underlying set $P$ can be endowed with the coordinatewise ordering, i.e. in 
$P\times P$ we have $$(p,q)\leq_{P\times P} (p',q') \textrm{ if and only if }
(p\leq_P p' \textrm{ and } q\leq_P q').$$
Since the ordering relation $\leq_P$ of $P$ is a subset of $P\times P$ 
it inherits the ordering described above such that it 
can be regarded as a poset in its own right. We denote this poset constructed
using $\leq_P$ as underlying set by $\Phi(P)$. Note that trivially by definition
$\Phi(P)$ is a subposet of $P\times P$ and we have $\Phi(P) = P\times P$ 
if and only if $P$ is a singleton.

If $P,Q$ are posets and $f:P\to Q$ is an order-preserving function, it is
easily seen that the restriction of $$(f\times f): P\times P \to Q\times Q
\textrm{ defined by } (p_1,p_2)\mapsto (f(p_1),f(p_2))$$ to $\Phi(P)$ gives
rise to an order-preserving function $$\Phi(f):\Phi(P)\to \Phi(Q).$$
It is easy to verify that with this construction we can make $\Phi$ into
a functor from the category of posets with order-preserving functions to itself.

Another easy calculation shows that if $L$ is a lattice then so is $\Phi(L)$.
Operations are componentwise; indeed $\Phi(L)$ is a sublattice
of the lattice $L\times L$. For $L\in \D$
it turns out that $\Phi(L)$ is a $(0,1)$-sublattice of $L\times L$ and
therefore $\Phi(L)\in \D$.

Perhaps not surprisingly, given a $(0,1)$-lattice homomorphism $f:L\to K$
between bounded (not necessarily distributive) lattices, the
map $\Phi(f): \Phi(L)
\to \Phi(K)$ is a lattice $(0,1)$-homomorphism as well. Routine verification
shows that $L \mapsto \Phi(L)$ and $f\mapsto \Phi(f)$ gives rise to a functor
$\Phi:\D \to \D$. This is what we want to have a closer look at in the following.
In section \ref{calcsection} we calculate $\Ip(\Phi(L))$ in terms of $\Ip(L)$
and in section \ref{priestleyPhi} we look at the interaction of $\Phi$ with
Priestley duality and 
characterise those bounded distributive lattices $L$ such that there is
$K$ with $\Phi(K) \cong L$.

A similar and in some way more general construction was studied by
J.D.~Farley in \cite{Jonathan}.

%------------------------------------
\section{Calculating $\Ip(\Phi(L))$}\label{calcsection}
In this section we express the collection of prime ideals of $\Phi(L)$ in 
terms of $\Ip(L)$.
\begin{lemma}\label{prodlem}
If $S$ is a ideal of $\Phi(L)$ then $$S=(pr_1(S)\times pr_2(S))
\cap \Phi(L)$$ where $pr_j: L\times L\to L$ is defined by $(l_1, l_2)
\mapsto l_j$ for $j=1,2$.
\end{lemma}
\begin{proof}
Certainly $S \subseteq (pr_1(S)\times pr_2(S))
\cap \Phi(L)$. On the other hand suppose that $(a,b) \in (pr_1(S)\times pr_2(S))\cap \Phi(L)$. So there is $b_a, a_b\in L$ such 
that $(a,b_a), (a_b, b)\in S$. Therefore $(a\vee a_b, b\vee b_a)\in S$ which entails $(a,b)\in S$, since $S$ is an ideal.
\end{proof}
For notational convenience, let $S_i$ denote $pr_i(S)$ for $i=1,2$. Note that $S_1=\{a\in L : (\exists b \in L):(a,b)\in S\}$. For $S_2$, a similar statement holds.
\begin{lemma}\label{primlem}
If $S$ is a prime ideal of $\Phi(L)$ then $S_1\in \Ip(L)$
and $S_2\in \Ip(L)\cup\{L\}$.\end{lemma}
\begin{proof}
First note that $S_1 \neq L$: for if we had $1\in S_1$, then there 
would be $b\in L$ such that $(1,b)\in S\subseteq \Phi(L)$, so $b=1$.
But if $S$ contains $(1,1)$ then we have $S=\Phi(L)$.

Moreover it is fairly easy to see that $S_1, S_2$ are ideals. Now suppose that $c,d \notin S_1$ but $c\wedge d \in S_1$. So there is
$b\in L$ such that $(c\wedge d, b)\in S$. Now $S\subseteq \Phi(L)$ entails $c\wedge d \leq b$. Since $S$ is a down-set of $\Phi(L)$, certainly $(c\wedge d, c\wedge d) \in S$. Moreover we have $(c,c), (d,d)\notin S$ (because $c,d \notin S_1$), so $(c,c)\wedge (d,d)= (c\wedge d, c\wedge d) \in S$ which contradicts $S$ being prime. With a similar argument we show that $S_2$ has the ''prime property'' - although it is possible that $S_2 = L$. 
\end{proof}
\begin{lemma}\label{eqlem} 
If $S$ is a prime ideal of $\Phi(L)$ then $S_2 
\in \{S_1, L\}$. 
\end{lemma}
\begin{proof}
First we show $S_2 \supseteq S_1$. Let $a\in S_1$, so by definition of $S_1$
there exists $b\in L$ such that $(a,b)\in S\subseteq \Phi(L)$. By construction
of $\Phi(P)$ this implies
$a\leq b$. Note that $S$ is a down-set of $\Phi(P)$ and $(a,a)\in \Phi(P)$
by definition of $\Phi$. Moreover, $(a,a)\leq (a,b)$ in $\Phi(P)$ since
$\Phi(P)$ is ordered coordinatewise. Because $S$ is a down-set, one obtains
$(a,a) \in S$ and therefore $a\in S_2$ by definition of $S_2$.

Now suppose that $S_2$ is a proper superset of $S_1$. We want to show that $S_2=L$. Suppose $1\notin S_2$. Take $y\in S_2 \setminus S_1$. There is $a\in L$ such that $(a,y)\in S$ (in particular $a\leq y$). So $(y,y)\notin S$ and $(a,1) \notin S$ (because $1\notin S_2$), but $(a,1)\wedge (y,y) =
(a,y)\in S$, contradicting $S$ being prime. 
\end{proof}
\begin{corollary} For the lattice $L$ we have,
$$\Ip(\Phi(L))=\{(I\times I)\cap \Phi(L); I\in \Ip(L)\} \cup
\{(I\times L) \cap \Phi(L); I\in \Ip(L)\}.$$
\end{corollary}
\begin{proof}
It is straightforward to check that $(I\times I)\cap \Phi(L)$ and
$(I\times L)\cap \Phi(L)$ are prime ideals of $\Phi(L)$ whenever $I$ is a prime ideal of $L$. 

On the other hand, suppose that $S\in \Ip(\Phi(L))$. By Lemma \ref{prodlem} the prime ideal $S$ can be written as $(S_1\times S_2)\cap \Phi(L)$. 
From Lemma \ref{primlem} we get that $I:=S_1$ is prime and finally
Lemma \ref{eqlem} implies that $S$ is either $(I\times I)\cap \Phi(L)$ or $(I\times L)\cap \Phi(L)$.
\end{proof}
\def\2{\underline{2}}
%-----------------------------------------------
\section{When is $L$ isomorphic to $\Phi(K)$ for some $K$?}\label{priestleyPhi}
The following question arises naturally: When is a bounded distributive lattice $L$ isomorphic to $\Phi(K)$ for some bounded
distributive lattice $K$? One special Priestley space will be the key here.
Denote by $\2$ the ordinal $2=\{0,1\}$ with its standard ordering and the
discrete topology.

With the aid of Priestley duality we are able to give an answer to that question. Let the pair of functors be denoted by $\X:\D\to {\mathcal P}$ and $\E:{\mathcal P} \to \D$ where ${\mathcal P}$ denotes the category of Priestley spaces with order-preserving continuous functions.
\begin{lemma}\label{priestleylem}
Let $X$ be a Priestley space. Then 
\begin{enumerate}
\item $\Phi(\E(X))\cong \E(X\times \2)$ in $\D$ and dually
\item $\X(\Phi(L)) \cong \X(L) \times \2$ in ${\mathcal P}$.
\end{enumerate}
\end{lemma}
\begin{proof} For the first statement, consider the function
$$\varphi:\Phi(\E(X))\to \E(X\times \2) \textrm{ defined by }
(d,e) \mapsto d\times \{1\} \cup e\times \{0\}$$
for clopen down-sets $d,e$ of $X$, and also
$$\psi: \E(X\times \2)\to \Phi(\E(X)) \textrm{ defined by }
c \mapsto (c_1, c_0)$$
for each clopen down-set $c$ of $X$ where $c_i:=\{x\in X: (x,i)\in c\}$ for $i=0,1$. We claim that
$\varphi$ and $\psi$ are order-preserving inverses of each other and
therefore provide an order (and lattice) isomorphism between 
$\Phi(\E(X))\cong \E(X\times \2)$. First note that for 
$(d,e) \in \Phi(\E(X))$ we have $d\subseteq e$ and therefore 
$\varphi((d,e))=(d\times \{1\}) \, \cup \, (e\times \{0\})$ is a clopen 
down-set in $X\times \2$. On the other hand, if $c$ is a clopen down-set 
of $X\times \2$ then $c_0 \supseteq c_1$, and clearly $c_0, c_1$ are clopen 
down-sets of $X$, so $(c_1, c_0) \in \Phi(\E(X))$. It is straightforward to 
check that $\varphi$ and $\psi$ are both order-preserving, so it remains to 
show that they are inverses of each other. Note that $\varphi(d,e)_1=d$ and 
$\varphi(d,e)_0=e$ for clopen down-sets $d,e$ of $X$. So 
$\psi(\varphi(d,e))=(d,e)$. Moreover for any clopen down-set $c$ of 
$X\times\2$ we have $c=(c_1\times\{1\})\, \cup \, (c_0\times \{0\})$, 
so $\varphi(\psi(c))=c$.

As for the second statement, let $X:=\X(L)$. If we apply the functor $\X$ to statement 1, we get $\X(\Phi(\E(X)))\cong \X(\E(X\times \2))$.
So we get with that and Priestley duality:
$$ \X(\Phi(L))\cong\X(\Phi(\E(X)))\cong \X(\E(X\times \2))\cong
X\times\2$$ which proves statement 2.
\end{proof}
\begin{theorem}\label{characterisation}For $L\in \D$ the following statements 
are equivalent:
\begin{enumerate}
\item $L$ is isomorphic to $\Phi(K)$ for some bounded
distributive lattice $K$
\item The Priestley space $\X(L)$ is order-homeomorphic to
$Y\times \2$ for some Priestley space $Y$. 
\end{enumerate}
\end{theorem}
\begin{proof}
Let $L\cong \Phi(K)$. Then by Lemma \ref{priestleylem}, statement 2, we get
$$\X(L)\cong\X(\Phi(K))\cong \X(K)\times \2.$$ So, taking
$Y:=\X(K)$ we are done.

For the other directition, suppose $\X(L)\cong Y\times \2$. By
Lemma \ref{priestleylem}, statement 1, we get
$$L\cong\E(\X(L))\cong\E(Y\times\2)\cong \Phi(\E(Y)).$$ So, taking
$K:=\E(Y)$ we are done.
\end{proof}
%-----------------------------------------------------
\section{Fixed points of $\Phi$}
Another natural question arising in the context of the functor
$\Phi$ is finding fixed points of $\Phi$, that is, distributive
$(0,1)$-lattices $L$ with the property that $\Phi(L) \cong L$. 

With Theorem \ref{characterisation} we can say
\begin{quote} A lattice $L$ is isomorphic to $\Phi(L)$ if and only if
for the Priestley space $Y$ of $L$ has the property that
$Y\cong Y\times \2$ (meaning there is an order-preserving 
homeomorphism from $Y$ to $Y\times \2$).
\end{quote}

The search for Priestley spaces $Y$ with the property $Y\cong Y\times \2$
gives rise to an example of a fixed point of $\Phi$.

\begin{example} For the Priestley space $Y=\2^\omega$ (endowed 
with the product topology and the coordinatewise ordering) we have
$Y\cong Y\times\2$.
\end{example}
The following map provides an order-preserving bijection from $Y$ to
$Y\times \2$: $$y \mapsto ( {\rm ls}(y),y(0)),$$ where ls denotes 
the left shift ${\rm ls}:\2^\omega\to \2^\omega$ given by ${\rm ls}(y)(n)=y(n+1)$
for all $n\in \omega$ and $y\in \2^\omega$. It is easy to see that the product topology
on $\2^\omega$ coincides with the interval topology on the poset $\2^\omega$
which is true for $\2^\omega \times \2$ as well. Recall that
the {\it interval topology} on any poset $P$ is the topology generated
by $$\{P\setminus [x,y]: x, y\in P \textrm{ and } x\leq y\}.$$ 
Recall that $[x,y] = \{z\in P: x\leq z\leq y\}$.
Note that
any order-isomorphism between posets is a homeomorphism between
the ground sets endowed with the interval topology. So the order-isomorphism
from above is a homeomorphism as well, which proves that
$Y$ and $Y\times \2$ are homeomorphic in $\P$. 

Applying Priestley duality to this example implies that for the lattice
$L = \E(Y)$ we have $\Phi(L)\cong L$. The object $L$ is
(isomorphic to) the free distributive
$(0,1)$-lattice generated by countably many points.

It is unclear how to characterise those Priestley spaces $Y$ with
$Y\cong Y\times \2$. The functor $\Phi$ and its fixed points gives rise
to more questions. 
\begin{question}
Does $\2^\omega$ embed into every Priestley space $Y$ having the property that
$Y\cong Y\times \2$? If not, is there a countable such Priestley space $Y$?
\end{question}
Of course, the functor $\Phi$ can be studied in the more general settings
of posets (even of preordered sets). Note that if $P$ and $Q$ are posets
which are fixed points of $\Phi$, then so is their disjoint union. So it
is more rewarding to consider connected posets only. Recall that a poset
$(P,\leq)$
is connected if $$tr(\leq \cup (\leq)^{-1})) = P\times P$$ where
$(\leq)^{-1} = \{(y,x): x\leq y\}$ and $tr$ denotes the transitive closure.
\begin{question} Is there a connected poset $P$ with more than one point 
such that 
$P$ is not a lattice and $\Phi(P)\cong P$?
\end{question}

There are several ``cardinal functions'' in the category of posets. Let
us just mention the order dimension and the width. The {\it width} is
the supremum of all cardinalities of anti-chains of a poset $(P,\leq)$ where
an anti-chain is a subset $A \subseteq P$ such that $x\neq y \in A$ 
implies $x\not \leq y$ and $y\not\leq x$. Moreover recall that any ordering
relation equals the intersection of all total ordering relations containing
it. (In a {\it total} ordering relation we have $x\leq y$ or $y\leq x$ 
for all $x,y$ in the ground set.) The {\em order dimension} 
of a poset $(P,\leq)$
is the minimal cardinality $\kappa$ such that 
there is a collection ${\mathcal S}$ of 
total ordering relations such that the intersection of ${\mathcal S}$ equals
the given ordering relation and ${\rm{card}}({\mathcal S}) = \kappa$.

Natural questions arise
when looking at those functions' interaction with $\Phi$, especially
in the case of finite posets. One example would be:
\begin{question}If $P$ is a finite poset, how does its order dimension 
compare to that of $\Phi(P)$?\end{question}

\section{Acknowledgement} We wish to thank Jonathan
Farley for helpful 
discussions about the functor $\Phi$ in the context of bounded distributive
lattices. Moreover, we are grateful to the (anonymous) referee for pointing
out a number of mathematical errors.

%------------------------ now paper ends ---------------------

\bibliographystyle{amsplain}

\end{document}